\documentclass[reqno,12pt]{amsart}
\usepackage{a4wide}
\usepackage{mathrsfs}

\numberwithin{equation}{section}

\newtheorem{theo}{Theorem}

\theoremstyle{remark}

\def\restriction#1#2{\mathchoice
	{\setbox1\hbox{${\displaystyle #1}_{\scriptstyle #2}$}
		\restrictionaux{#1}{#2}}
	{\setbox1\hbox{${\textstyle #1}_{\scriptstyle #2}$}
		\restrictionaux{#1}{#2}}
	{\setbox1\hbox{${\scriptstyle #1}_{\scriptscriptstyle #2}$}
		\restrictionaux{#1}{#2}}
	{\setbox1\hbox{${\scriptscriptstyle #1}_{\scriptscriptstyle #2}$}
		\restrictionaux{#1}{#2}}}
\def\restrictionaux#1#2{{#1\,\smash{\vrule height .8\ht1 depth .85\dp1}}_{\,#2}}

\def\cC{\mathcal{C}}
\def\de{\delta}

\def\ze{\zeta}

\def\Ga{\Gamma}

\def\({\left(}
\def\){\right)}
\def\[{\left[}
\def\]{\right]}
\def\ii{\infty}

\newcommand\bin[2]{\left(\!\!\!\begin{array}{c}#1\\#2\end{array}\!\!\!\right)}

\def\dis{\displaystyle}

\def\N{\mathbb{N}}
\def\Q{\mathbb{Q}}
\def\Z{\mathbb{Z}}
\def\C{\mathbb{C}}

\def\R{\mathbb{R}}
\date{}

\begin{document}
\title{Remainder Pad\'e approximants   for  the Hurwitz zeta function}
\author{Marc Pr\'evost } 
\address{LMPA Joseph Liouville, Centre Universitaire de la Mi-Voix,
B\^at H. Poincar\'e,
50 rue F. Buisson,
BP 699,
62228 Calais Cedex}
\email{marc.prevost@univ-littoral.fr}

 

\begin{abstract}
Following our earlier research, we use the method  introduced by the author in \cite{prevost1996} named  Remainder Pad\'e Approximant in \cite{rivoalprevost}, to construct  approximations of the Hurwitz zeta function. We prove that these approximations are  convergent on the positive real line.
Applications to new rational approximations of $\ze(2)$ and $\ze(3)$ are given.

 \end{abstract} 

\maketitle

\section{Introduction}

In \cite{prevost1996}, we gave a new proof of the irrationality of $\zeta(2)=\sum_{k=1}^{\infty}1/k^2$ (and also of $\zeta(3)=\sum_{k=1}^{\infty}1/k^3$) based on an explicit construction of some Pad\'e approximants of the remainder term
$R_2(1/n)=\sum_{k=n}^{\infty}1/k^2$. Pore precisely, we have that $\zeta(2)=\sum_{k=1}^{n-1}1/k^2+R_2(1/n)$ with
$$R_2(z)=\sum_{k=0}^{\infty}\dfrac{z^2}{(z k+1)^2}.$$
The function $R_2$ is meromorphic in $\C\setminus \{0,-1,-\dfrac{1}{2},-\dfrac{1}{3},\cdots\}$ and so is not holomorphic at $0$.

However, it is $\cC^\infty$ at $z=0$ and admits a Taylor series $\hat R_2(z)=\sum_{k=0}^{\infty}B_k z^{k+1}$ which is an asymptotic development convergent only at $z=0$, (the radius of convergence is $0$), where  $B_k$ is the $k$th Bernoulli number.

The idea was to compute an explicit representation of diagonal Pad\'e approximant $P_n(z)/Q_n(z)=[n/n]_ {\hat R_2}(z)
\in \Q(z)$ of the series $\hat R_2(z)$ with a good estimation of the error term $E_n(z)=R_2(z)-[n/n]_ {\hat R_2}(z).$

At last, we find that
$$Q_n(1/n)E_n(1/n)=Q_n(1/n)\ze(2)-Q_n(1/n)\sum_{k=1}^{n-1}1/k^2-P_n(1/n)=q_n\ze(2)-p_n$$
provides a sequence of rational approximation $P_n/q_n$ which proves tje irrationlaity of $\ze(2).$
The surprise was to find that it is exactly the sequences used by Ap\'ery \cite{apery} for the same  purpose.

The same method applied to the series $\zeta(3)=\sum_{k=1}^{\infty}1/k^3$ provides also the Ap\'ery numbers 
given in  \cite{apery}.

The name  Remainder Pad\'e Approximant (RPA) has been used later in a paper written with T. Rivoal \cite{rivoalprevost} for the exponential function.

In this paper, we apply this method and  prove that some sequence of  RPA for Hurwitz zeta function is convergent on the real line.

For $a\in \C, \Re(a)>0$ , the   Hurwitz zeta function is defined as 

\begin{eqnarray}
\ze (s,a)& = &\dis{\sum ^\ii _{k=0}} \dis{\frac{1}{(k+a)^s}}
\label{1}
\end{eqnarray} 

\noindent where the Dirichlet series on the right hand side   of (\ref{1}) is
convergent for $\Re(s) > 1$ and uniformly convergent in any finite region  
where $\Re(s) \geq 1 + \de $, with $\de  > 0$. It defines an
analytic function for $\Re(s) > 1$. The Riemann $\ze$ function is $\ze(s,1)$.
   
An integral representation of $\ze(s,a)$ is
\begin{eqnarray}
\ze (s,a) = \dis{\frac{1}{\Ga(s)}} \int ^\infty _0
\displaystyle{\frac{x^{s-1}e^{-a x}}{1-e^{-x}}}dx \quad \Re(s) > 1,\Re(a)>0,
\end{eqnarray}

{where} \; $\Ga (s) =  \int ^\ii_0 y^{s-1} e^{-y}dy$ 
is the gamma function.

If $\Re(s) \leq 1$ the Hurwitz zeta function can be defined by the equation
:
$$\ze (s,a)=\Ga(1-s)\frac{1}{2 \pi i} \int_\cC \frac{z^{s-1}e^{a z}}{1-e^z}dz$$

where ${\cC}$ is some path in $\C$, which provides the
analytic continuation of $\ze (s,a)$ over the whole $s$-plane.

If we write the formula (\ref{1}) as :

\begin{eqnarray}
\ze (s,a) = \dis{\sum ^{n-1}_{k=0} \frac{1}{(k+a)^s}}+
\dis{\sum ^\ii _{k=0}} \dis{\frac{1}{(n+a+k)^s}}
 \end{eqnarray}

 and set
 $$\Psi(s,x) : =  \dis t^{1-s}\sum ^\ii_{k=0}\({\frac{t}{1+kt}}\)^s$$

then
\begin{eqnarray}
 \ze (s,a) = \dis{\sum ^{n-1} _{k=0}}
\dis{\frac{1}{(k+a)^s}} + (n+a)^{1-s}
\Psi\(s,\frac{1}{n+a}\).
\label{formzeta}
 \end{eqnarray}

Another expression of $\Psi(s,t)$ is
 \cite[chapt. XIII]{watson-witt}

\begin{eqnarray}
\Psi(s,t) = \dis{\frac{t^{1-s}}{\Ga (s)}}\int ^\ii _0 u^{s-1}\dis{
\frac{e^{-u/t}}{1-e^{-u}}}du
\end{eqnarray}

\noindent whose Taylor series is

\begin{eqnarray}
\Psi(s,t) = \dis{\sum ^\ii_{k=0}}
\dis{\frac{B_k}{k!}} (s)_{ k-1} (-1)^k t^{k},\label{DAS}
 \end{eqnarray} convergent only at $t=0$,
  where $B_k$ is the $k$th   Bernoulli number and $(s)_k$ is the Pochammer symbol
($(s)_p:=s(s+1)(s+2)\cdots(s+p-1),(s)_{-1}=1/(s-1)$).

We replace in the remainder term in  (\ref{formzeta}), $\dis\Psi\(s,\frac{1}{n+a}\)$
by some   Pad\'e approximant $[m_1/m_2]_{\Psi(s,t)}=\widetilde
 Q_{m_1}(t)/\widetilde
 P_{m_2}(t)$ to the function $\Psi(s,t)
$ computed at $ t={1}/(n+a)$.  We obtain the RPA approximant
\begin{eqnarray}
  RPA{(n,m_1,m_2)}:={\sum ^{n-1} _{k=0}}
{\frac{1}{(k+a)^s}} + (n+a)^{1-s}\dfrac{\widetilde Q_{m_1}({1}/{(n+a)})}{\widetilde P_{m_2}({1}/{(n+a)})}.
\end{eqnarray}

 The crucial point is the convergence of this sequence $RPA{(n,m_1,m_2)}$ when one or more of the parameter
 $n,m_1,m_2$ tends to infinity. \\The convergence of  Pad\'e approximants is proved for many functions as
 meromorphic function, Stieltjes function, etc...(\cite{bgm1996}).
 In this paper, we   prove that the function $\Psi\(s,t\)$ is a Stieltjes function in the variable $t$ when $s$ is a real positive number.
 First we have to recall the  definition and some properties of Pad\'e approximants

 \section{Pad\'e approximants}
 
 Let $f$ be a function whose Taylor expansion about $t = 0$ is
 
 $$f(t)=\sum ^\infty_{i=0} c_i  t^i.$$
 
 The Pad\'e approximant $[m_1/m_2]_f$ to $f$ is a rational fraction $\widetilde
 Q_{m_1}/\widetilde P_{m_2}$ whose denominator has degree $m_2$ and whose
 numerator has degree $m_1$ so that its expansion in ascending powers
 of $t$ coincides with the expansion  of $f$ up to the
 degree $m_1+m_2$ , i.e.
 
 \begin{eqnarray}
 \widetilde P_{m_2}(t) f(t) - \widetilde Q_{m_1}(t) = 0(t^{m_1+m_2+1})\;,\quad t
 \rightarrow 0 \end{eqnarray}

 The theory of Pad\'e approximation is linked with the theory of orthogonal
 polynomials as following \cite{PTA}: let us define the linear functional $c$ acting on ${\mathcal P}$,
 space of polynomials by 
 
 \begin{eqnarray}
 &c : &{\mathcal P}\rightarrow \R  (or \; \C )\\
 &&x^i \rightarrow <c, x^i> = c_i \quad i = 0, 1, 2, \ldots
 \end{eqnarray}
 
 and
 
 \begin{eqnarray}
 \mathrm{if} \; p \in {\Z}\;\;  <c^{(p)}, x^i> := <c,x^{i+p}> = c_{i+p} \quad i = 0, 1, 2,
 \ldots (c_i=0, i<0)
 \end{eqnarray} 
 then the denominators of the
 Pad\'e approximants $[m_1/m_2]_f$ satisfy the following orthogonality
 property :
 
 $$<c^{(m_1-m_2+1)}, x^i P_{m_2}(x)> = 0 \quad i = 0,1,2,\ldots,n-1$$
 \noindent where $P_{m_2}(x) = x^{m_2} \widetilde P_{m_2}(x^{-1})$.
 
 Let us  define the associated polynomials :
 
 \begin{eqnarray}
 R_{m_2-1}(t) := <c^{(m_1-m_2+1)},\frac{P_{m_2}(x)-P_{m_2}(t)}{x-t}>, R_{{m_2}-1}\in {\mathcal P}_{{m_2}-1}\label{associe}
 \end{eqnarray}
 \noindent where $c^{(m_1-{m_2}+1)}$ acts on the variable $x$ . Then 
 
 \begin{eqnarray}
 \widetilde Q_{m_1}(t) =\( {\sum ^{{m_2}-m_1}_{i=0}} c_it^i\)
 \widetilde P_{m_2}(t) + t^{m_1-{m_2}+1}\widetilde R_{{m_2}-1}(t)
 \end{eqnarray}
 \noindent where $\widetilde R_{{m_2}-1}(t) = t^{{m_2}-1} R_{{m_2}-1}(t^{-1})$ and $c_j = 0$
 for $j < 0.$
 
 If $c$ admits an integral representation with a function $\alpha $ non
 decreasing, with bounded variation, 
 
 \begin{eqnarray}
 c_i = \int _{L} x^i d\alpha (x),\label{integralrepr}
 \end{eqnarray}
 \noindent then   the function $f$ is a Stieltjes function and the theory of Gaussian quadrature shows that the polynomials
 $P_n$ orthogonal with respect to $c$, have all their roots in the
 support of the function $\alpha $.
Moreover, the convergence is proved if the coefficients $c_i$ satisfy the Carleman condition \cite[p.240]{bgm1996}.

The aim of the paper is to find, in the case of Hurwitz zeta function,  the  weight function  $d\alpha$ depending on $s$ and prove that it is a positive  function. 

The error is defined by
 \begin{eqnarray}
 \mbox{\underline{error}}\; : f(t) - [m_1/m_2]_f(t) =
 \displaystyle{\frac{t^{m_1+m_2+1}}{\widetilde P^2_{m_2}(t)}}c^{(m_1-{m_2}+1)} \left (
 \displaystyle{\frac{P^2_{m_2}(x)}{1-xt}}\right ).
 \end{eqnarray}

 The above expression of the error is understood as a formal one if $c$ is only
 a formal linear functional \cite[chapt. 3]{PTA},  but if $c$ admits the
integral representation (\ref{integralrepr}) then the error becomes :
 
 \begin{eqnarray}
 f(t) - [m_1/m_2]_f(t) = f(t) - {\frac{\widetilde Q_{m_1}(t)}{\widetilde P_{m_2}(t)}} = {\frac{t^{m_1+{m_2}+1}}{\widetilde P^2_{m_2}(t)}} 
 \int _{L} 
 x^{m_1-m_2+1} \frac{P^2_{m_2}(x)}{1-x t}d\alpha (x)
 \end{eqnarray}
 
 In the particular case $m_1 = m_2-1$,
 
 \begin{eqnarray}
 f(t) - [m_2-1/m_2]_f(t) = \displaystyle{\frac{t^{2m_2}}{\widetilde
 		P^2_{m_2}(t)}}\int _{L}
 {\frac{P^2_{m_2}(x)}{1-xt}} d\alpha (x).\label{errorpade}
 \end{eqnarray}
 
 Note that 
 
 $$[{m_2}/{m_2}]_f(t) = c_0+t [{m_2}-1/{m_2}]_{f_1}$$
 
 \noindent and  
 
 $$[m_2+p/m_2]_f(t) = c_0+c_1 t+\cdots +c_{p} t^{p}+t^{p+1} [m_2-1/m_2]_{f_p}(t)$$
 
 where $$f_p(t)=\sum ^\infty_{i=0} c_{i+p+1 }  t^i.$$

 The orthogonal  polynomial can be expressed with determinant:
 
 Let $(e_k)_{k\in \N}$ a basis of the space of polynomials. If 
 $$\left\langle c,e_i\;e_j\right\rangle =\alpha_{i,j}$$
 then
 \begin{eqnarray}
 P_n(x)&=&k_n\left|  \begin{array}{cccc}
\alpha_{0,0}&\alpha_{1,0}&\cdots&\alpha_{n,0}\\
\alpha_{0,1}&\alpha_{1,1}&\cdots&\alpha_{n,1}\\
\hdots&\hdots&\hdots&\hdots\\
\alpha_{0,n-1}&\alpha_{1,n-1}&\cdots&\alpha_{n,n-1}\\
e_0(x)&e_1(x)&\cdots&e_n(x)
 \end{array}\right|  \label{hadamard}
 \end{eqnarray}
 where $k_n$ is some constant.
 
 This property will be used in section \ref{particular case}.

 \section{Statements of the results}
 In this section, we will prove that it is possible to compute the weight function underlying the coefficients 
 $\frac{B_k}{k!} (s)_{ k-1} (-1)^k$
of the function 
 $\Psi(s,t) = {\sum ^\ii_{k=0}}
 {\frac{B_k}{k!}} (s)_{ k-1} (-1)^k t^{k}.$
 The weight function $w_s$  depends on $s$. Thus, for particular values of $s$,  $\Psi(s,t)$ is a Stieltjes function 
 in the variable $t$ and the convergence of RPA will be proved for  $t=1/(n+a)$.

\begin{theo}

For any complex numbers $ s,a $, $\Re(s)>-1,\Re(a)> 0$  and any positive integer $m$, $m>\Re(s)-1$, then

\begin{eqnarray}
\zeta(s,a)=\frac{1}{ a^{s-1}}\left(\frac{1}{s-1}+\frac{1}{2 a}+
\int_0^\infty
\frac{  1}{a^2+ x^2} w_s(x)dx\right)\label{zeta}
\end{eqnarray}

where the weight $w_s$ is defined by:

$$
 w_s(x):=\frac{2 (-1)^{m}x^s}{\Gamma(s)\Gamma(m+1-s)} \int_x^\infty (t-x)^{m-s}\frac{d^m}{dt^m}\left(\frac{1}
{e^{2\pi t}-1}\right)dt.
$$

\label{theo1}
\end{theo}

We will prove in the section \ref{section3}   that this  formula (\ref{zeta}) is a consequence of Hermite's formula for the function $\ze(s,a)$.

As found by Touchard \cite{touch}, Bernoulli numbers satisfy 
$$B_k=-i \frac{\pi}{2}\int_L x^k \frac{dx}{\sin^2(\pi x)},$$
where $L$ is the line  $L:=-\frac{1}{2}+i \R$
and thus can be viewed  as moments of order $k\geq 0$  of  the  positive weight function $1/{\sin^2(\pi x)} $ on the line $L:=-\frac{1}{2}+i \R$:
So  for $s=2$, the coefficients ${\frac{B_k}{k!}} (s)_{ k-1} =B_k$ of  (\ref{DAS}) are moments of a positive weight.
For $s=3$, these coefficients appear as their derivative $(k+1)B_k$. The derivative of $1/{\sin^2(\pi x)}$ equals
is a weight function symmetric around $-1/2$ whose support is also the line $L$.
But for integer value of $s$ greater than 4, it is no more true for all $k\geq 0$.
A much more difficult case is the case $s$ is real.

From the previous Theorem, we obtain an integral representation of $\dis{\frac{B_k}{k!}} (s)_{ k-1} $ but for $k\geq 2$.

\begin{theo}Integral representation of Bernoulli numbers.\label{theo2}

If $s$ a complex number  such that $\Re(s)>-1$, then

\begin{eqnarray}
B_{k+2}\frac{(s)_{k+1}}{(k+2)!}&=&(-1)^{k/2}
\int_{0}^{+\infty}
x^{k} w_s(x)dx \hskip 1 cm (k \;even)  \label{bernoulli}\\
&=& 0  \hskip 5 cm\;(k\; odd)
\end{eqnarray}

 with 

\begin{equation} \dis w_s(x):=\frac{2 (-1)^{m}x^s}{\Gamma(s)\Gamma(m+1-s)} \int_x^\infty (t-x)^{m-s}\frac{d^m}{dt^m}\left(\frac{1}{e^{2\pi t}-1}\right)dt
\end{equation}

{\rm where\; $m$ is\; any\; integer \:satisfying\;} $ m > \Re(s)-1.$

\label{teo2}
\end{theo}

\textbf{Proof}
We use the expansion of $\zeta(s,a)$ in terms of Bernoulli numbers \cite [p.160]{srivastava-choi}.
$$\zeta(s,a)=a^{-s}+\sum_{k=0}^n (s)_{k-1}\frac{B_k}{k!}a^{-k-s+1}+\frac{1}{\Gamma(s)}\int_0^\infty \left(\frac{1}{e^t-1}-\sum_{k=0}^n \frac{B_k}{k!}t^{k-1}\right)e^{-a\;t}\;t^{s-1}\;dt$$
valid for $\Re(s)>-2n+1,\Re(a)>0$ and $n$ a positive integer.

By identification with formula (\ref{zeta}), we get the formula (\ref{bernoulli}).

{\bf Remark}: The positivity of the weight function is important because it will imply the convergence of Pad\'e approximants.  In the case where $s$ is a positive real number then $w_s$ is positive on its support. This gives the   following main theorem.

 \begin{theo}
 	
 	For all positive real number $s$, for all complex number $a, \Re(a)>0$, for all integers $n\geq 0, p\geq 1$, the following sequence    
 	
 		\begin{eqnarray}RPA{(n,m+p,m)}:={\sum ^{n-1} _{k=0}}
 	{\frac{1}{(k+a)^s}}  +
 	(n+a)^{1-s}\restriction{[m+p/m]_{\Psi(s,t)}}{t=1/(n+a)}  .\label{RPA approx}
 	\end{eqnarray}
 	converges to $\ze(s,a)$ when $m$ tend to infinity.\label{teo4}
 \end{theo}

 \vskip 10pt
  These Theorems will be proved in Section \ref{section3}. 
  
  \vskip 10pt
{\bf Remark}

For $n=0,p=-1$, this result is proved  only for $s=2, 3$ in \cite{prevost1996}.

If $n=0$, then formula (\ref{RPA approx}) reduces to

$$\forall s >0 ,\forall a \in \C, \Re(a)>0, \forall p\geq 1,  \lim_{m\rightarrow \infty} a^{1-s}\restriction{[m+p/m]_{\Psi(s,t)}}{t=1/a}=\ze(s,a).$$

Moreover, for the particular case $n=0$, $a=1$, then $RPA{(0,m+p,m)}$ is a rational function depending on the variable $s$ and  which converges to $\ze(s,1)=\ze(s)$ when $m$ tends to $\ii$ (by Theorem \ref{teo4}).

For example, for $p=1$, 
\\$RPA(0,2,1)=\frac{(s+2) (s+3)}{12 (s-1)} ,\\RPA(0,3,2)=\frac{(s+2) \left(s^2+12 s+31\right)}{2 (s-1) \left(s^2+3 s+62\right)},\\RPA(0,4,3)=\frac{(s+2) (s+4) \left(3 s^2+88 s+405\right)}{120 (s-1) \left(s^2+7 s+54\right)},\\RPA(0,5,4)=\frac{(s+2) (s+4) \left(s^5+124 s^4+2644 s^3+23730 s^2+92939 s+122130\right)}{2 (s-1) \left(s^6+101 s^5+1755 s^4+21415 s^3+153884 s^2+657564 s+977040\right)} \\
RPA(0,6,5)=\frac{(s+2) (s+4) (s+6) \left(s^5+746 s^4+28162 s^3+498112 s^2+3612925 s+8457750\right)}{84 (s-1) \left(s^6+429 s^5+10387 s^4+134511 s^3+1044772 s^2+4891020 s+9666000\right)}$

\vskip 1cm

 If $n\geq 1$ and $a=1$, then $RPA(n,m+p,m)$ is no more a rational fraction since it  contains powers of $(n+a)$ with $s$ as exponent:
\\$RPA(1,2,1)=\frac{2^{1-s} \left(s^2+11 s+36\right)}{48 (s-1)}+1,\\
RPA(1,3,2)=\frac{2^{-s-1} \left(s^3+26 s^2+231 s+726\right)}{(s-1) \left(s^2+3 s+242\right)}+1,\\
RPA(1,4,3)=\frac{2^{-s-4} \left(3 s^4+166 s^3+2937 s^2+22334 s+64800\right)}{15 (s-1) \left(s^2+7 s+180\right)}+1,
\\
RPA(2,2,1)=\frac{3^{1-s} \left(s^2+17 s+90\right)}{108 (s-1)}+2^{-s}+1,\\
RPA(2,3,2)=\frac{3^{-s} \left(s^3+38 s^2+527 s+2710\right)}{2 (s-1) \left(s^2+3 s+542\right)}+2^{-s}+1.$

\vskip 10pt
Nevertheless, if $s$ is an integer, if $a=1$ and $p\geq 1$, we obtain a sequence of rational approximations of $\ze(s)$ by computing
$RPA(n,m+p,m)$, depending on the two integers parameters $n$ and $m$.

\vskip 10pt

\section{Proof of Theorems}\label{section3}

\subsection{Proof of Theorem \ref{theo1}}

We start from   Hermite's formula for $\zeta(s,a),\Re(a)>0,$ which is a consequence of  Plana's summation formula:
\begin{equation}
\zeta(s,a)=\frac{1}{2}\;a^{-s}+\frac{a^{1-s}}{s-1}+2\int_0^\infty
(a^2+y^2)^{-s/2}\sin\left(s \arctan{\frac{y}{a}}\right)\frac{dy    }{e^{2\pi y}-1}.
\label{hermite}
\end{equation}
Note that this integral converges for all complex number $s\neq1$.

Let us define
$$I_s:= \int _0^t\frac{x^s(t-x)^{-s}}{1+x^2}dx \;\;\;(t\geq 0,-1<\Re(s)<1).$$

Using  the identity $$\int_0^\infty\frac{v^{s-1}}{1+\alpha\; v}={\alpha^{-s}}{\Gamma(s)\Gamma(1-s)}$$
we get, with the  change of variable $x=\dis\frac{t\;v}{1+v}$,

\begin{eqnarray*}
I_s&=&t\int_0^\infty \frac{v^s}{(1+v)^2+t^2v^2}dv\\&=&\frac{1}{2i}\int_0^\infty v^{s-1}\left( \frac{1    }{1+v(1-i t)}-
\frac{1    }{1+v(1+i t)}\right)dv\\&=&
\Gamma(s)\Gamma(1-s)
\frac{1}{2i}\left((1-i t)^{-s}-(1+i t)^{-s}\right).
\end{eqnarray*}

Setting $t=\tan \theta $, we obtain
$$I_s=\Gamma(s)\Gamma(1-s)\cos^s( \theta) \sin(s \theta)=\Gamma(s)\Gamma(1-s)(1+t^2)^{-s/2}\sin (s \arctan t).$$

Thus, the following identity
$$\frac{1}{\Gamma(s)\Gamma(1-s)}\int _0^t\frac{x^s(t-x)^{-s}}{1+x^2}dx=(1+t^2)^{-s/2}\sin (s \arctan t)$$

holds for all real $t\geq 0$ and complex $s$, $-1<\Re(s)<1$.

Now, for $s$ complex satisfying $\Re(s)>-1$ and $m$ an integer such that $m>\Re(s)-1$, by recurrence on $m$, it is easy to prove the following formula,

\begin{equation}
\frac{\sin(s \arctan t)}{(1+t^2)^{s/2}}=\frac{1}{\Gamma(s)\Gamma(m+1-s)}\frac{d^m}{dt^m}\left(\int_0^t
\frac{x^s(t-x)^{m-s}    }{1+x^2}\;dx\right).
\label{sin}
\end{equation}

To prove  Theorem \ref{theo1}, we replace $\displaystyle \frac{\sin(s \arctan (t/a))}{(a^2+t^2)^{s/2}}$ in Hermite's formula (\ref{hermite}), then apply $m$ integration by parts and permutation of the integrals.

\begin{eqnarray*}
\zeta(s,a)&=&\!\!\!\frac{a^{-s}}{2}\!+\!\frac{a^{1-s}}{s\!-\!1}\!+\!\frac{2\;a^{1-s}}{\Gamma(s)\Gamma(m\!+\!1\!-\!s)}\int_0^\infty
\frac{d^m}{dt^m}\left(\int_0^t
\frac{x^s(t\!-\!x)^{m-s}}{a^2+x^2}dx\right)\frac{dt    }{e^{2\pi t}\!-\!1}\\
&=&\!\!\!\frac{a^{-s}}{2}\!+\!\frac{a^{1-s}}{s\!-\!1}\!+\!\frac{2(-1)^ma^{1-s}}{\Gamma(s)\Gamma(m\!+\!1\!-s\!)}\int_0^\infty\!\!\!
\left(\int_0^t
\!\!\frac{x^s(t\!-\!x)^{m-s}}{a^2\!+\!x^2}dx\right)\!\frac{d^m}{dt^m}\!\!\left(\!\!\frac{1 }{e^{2\pi t}\!-\!1}\!\!\right)dt\\
&=&\!\!\!
\frac{a^{-s}}{2}\!+\!\frac{a^{1-s}}{s-1}\!+\!\frac{2(-1)^ma^{1-s}}{\Gamma(s)\Gamma(m\!+\!1\!-\!s)}\int_0^\infty\!\!
\!\!\!\frac{x^sdx}{a^2\!+\!x^2}
\int_x^\infty\!\!(t\!-\!x)^{m-s}
\frac{d^m}{dt^m}\!\!\!\left(\!\!\frac{1    }{e^{2\pi t}\!-\!1}\!\!\right)dt.\label{zeta(s,a)}
\end{eqnarray*}

The permutation in the last integral is valid since the function $
\frac{x^s}{a^2+x^2}
 (t-x)^{m-s}
\frac{d^m}{dt^m}\left(\frac{1    }{e^{2\pi t}-1}\right)$
is less than $ \frac{t^m}{a^2}\frac{d^m}{dt^m}\left(\frac{1    }{e^{2\pi t}-1}\right)$ which integrable on
$[0,\infty[$.

\vskip 10pt

\textbf{Remark}

If we consider  $a=1$, we get an integral representation of the Riemann zeta function:

\begin{eqnarray}
\zeta(s)&=&\frac{1}{2}+\frac{1}{s-1}+\int_0^\infty\frac{ 1}{1+ x^2} w_s(x)dx
\end{eqnarray}
with 

\begin{eqnarray}
\dis w_s(x):=\frac{2 (-1)^{m}x^s}{\Gamma(s)\Gamma(m+1-s)} \int_x^\infty (t-x)^{m -s}\frac{d^m}{dt^m}\left(\frac{1}{e^{2\pi t}-1}\right)dt
\end{eqnarray}

where $m$ is an  integer  satisfying  $ m > \Re(s)-1>0.$

So, if $s$ is an integer, we can take $m=s$, and

\begin{eqnarray}
\zeta(s)&=&\frac{1}{2}+\frac{1}{(s-1)}+ \frac{2 (-1)^{s-1}}{\Gamma(s)}
\int_0^\infty
\frac{  x^s}{1+ x^2} \frac{d^{s-1}}{dx^{s-1}}\left(\frac{1}{e^{2\pi x}-1}\right)dx\\
&=&\frac{1}{2}+\frac{1}{(s-1)}+ \frac{2}{\Gamma(s)}
\int_0^\infty
\frac{ 1}{e^{2\pi x}-1}\; \frac{d^{s-1}}{dx^{s-1}}\left(\frac{x^s}{1+x^2}\right)dx
\label{eq:}
\end{eqnarray}

\begin{eqnarray}
\zeta(2)&=&\frac{3}{2} +\pi
\int_0^\infty
\frac{  x^2}{1+ x^2}  \frac{1}{\sinh ^2 \pi x}dx= \frac{3}{2}+4\int_0^\infty
\frac{1}{e^{2\pi x}-1}\;\frac{x}{(1+x^2)^2}dx\\
\zeta(3)&=& 1+ \pi^2
\int_0^\infty
\frac{  x^3}{1+ x^2} \frac{\cosh \pi x}{\sinh ^3\pi x}dx=
1-2\int_0^\infty \frac{1}{e^{2\pi y}-1}\;\frac{x(x^2-3)}{(1+x^2)^3}dx\\
\zeta(4)&=&\frac{5}{6} + \frac{\pi^3}{3}
\int_0^\infty
\frac{  x^4}{1+ x^2}\frac{2+\cosh(2\pi x)}{\sinh^4\pi x}dx=\frac{5}{6}-8 \int_0^\infty \frac{1}{e^{2\pi x}-1}\;\frac{x(x^2-1)}{(1+x^2)^4}dx.
\end{eqnarray}

\subsection{Proof of Theorem \ref{theo3}}
We consider the function $\Psi(s,t)$ given in (\ref{DAS}) written as
\begin{eqnarray}
\Psi(s,t) =\frac{B_0}{s-1}-B_1\, t+t^{2}\Psi_2(s,t)
\end{eqnarray}
where $$\Psi_2(s,t):={\sum ^\ii_{k=0}}
{\frac{B_{k+2}}{(k+
	2)!}} (s)_{ k+1} (-1)^k t^k.$$

Thus the Pad\'e approximant $[m+p/m]_{\Psi(s,t)}$ satisfies
$$[m+p/m]_{\Psi(s,t)}=\frac{B_0}{s-1}-B_1 \,t+t^{2}[m+p-2,m]_{\Psi_2(s,t)}.$$

For $s$ positive real number, the weight $w_s$ is positive (Theorem \ref{theo2}) and so $\Psi_2(s,t)$ is a Stieltjes function.
Its coefficients $c_k:=\dis{\frac{B_{k+2}}{(k+1)!}} (s)_{ k+1} (-1)^k$ are positive and satisfy the Carleman condition, i.e. the series
$$\sum_{k=0}^{\ii}c_k^{-1/2k} \mathrm{diverges}.$$
Actually, the Bernoulli numbers of even  index satisfy (Bernoulli of odd index greater than 2 are zero)
$$\left| B_{2k}\right| \sim 4 \sqrt{\pi k}\(\frac{k}{\pi e}\)^{2k}$$
and so $$ c_{2k}^{-1/4k}\sim \(\frac{k}{\pi e}\)^{-1/2}$$ 
which is the general term of a divergent series.

We now apply the Theorem of \cite[p.240]{bgm} to conclude that for all integer $p, p\geq1,$
$$\lim_{m\rightarrow\ii}t^{s-1}[m+p/m]_{\Psi(s,t)}=t^{s-1}\Psi(s,t)=\ze(s,1/t)\mbox{ for all complex number  } t, \Re(t)>0$$

\section{Particular case: $s$ is an integer}\label{particular case}

If $s$ is an integer, we can improve the approximation of $\ze(s,a)$.
In (\ref{formzeta}), we replace $\Psi\(s,\frac{1}{n+a}\)$ by some Pad\'e approximant $[m+p/m]_{\Psi(s,t)}$ to obtain an approximation of $\ze (s,a)$:
 
$$ \ze (s,a) = \dis{\sum ^{n-1} _{k=0}}
 \dis{\frac{1}{(k+a)^s}} + (n+a)^{1-s}\lim_{m \rightarrow \infty}
 \restriction{[m+p/m]_{\Psi(s,t)}}{t=1/(n+a)}  .$$

For $s=2$ and $s=3$, we will find in this section the formal expression of these approximants and the expression of the error. 

  The case for $s\geq 4$ remains an open problem.

\subsection{Case $s=2$}

For $s=2$, the weight in the expression  (\ref{zeta}) is
\begin{eqnarray}
w_2(x):=2 x^2  \int_x^\infty  \frac{d^2}{dt^2}\left(\frac{1}
{e^{2\pi t}-1}\right)dt=\frac{\pi x^2  }{\sinh^2 (\pi x)}.
\end{eqnarray}
Another expression of $w_2(x)$ is 
$$w_2(x)=\frac{1}{\pi}\left| \Gamma(1+i x)\Gamma(1-i x)\right|  ^2.$$  
Generalizing a result by Carlitz \cite{carlitz}, Askey and Wilson \cite{askeywilson} gave an explicit expression:  
the orthogonal polynomial $P_n$ with respect to the weight function $w_2$ satisfies:
\begin{eqnarray}
P_m(x)&=&(m+1)(m+2)_3F_2\(\begin{array}{ccccc}
-m&m+3&1-x&\\
&&&;1\\
&2&2&
\end{array}\)\\&=&
\sum_{k=0}^{m}\(\begin{array}{c}
m+1\\
k+1
\end{array}\)
\(\begin{array}{c}
m+k+2\\
k+1
\end{array}\)
\(\begin{array}{c}
x-1\\
k
\end{array}\),\label{polorths=2}
\end{eqnarray}

and

\begin{eqnarray}
\int_{i \R}P_n(x)P_m(x)\frac{\pi x^2}{\sin^2\pi x}dx&=&0 \;\;(n\neq m),\\
&=&\frac{(-1)^n 2(n+1)(n+2)}{2n+3}\;\;(n=m).
\end{eqnarray}

The roots of the polynomials $P_n$ are located on the imaginary axis since  the weight $\dis\frac{\pi x^2}{\sin^2\pi x}$ is positive on this line (see \ref{integralrepr} ).

The three terms recurrence relation is
\begin{eqnarray}
P_{m+1}(x)&=&\frac{2(2m+3)}{(m+1)(m+2)}\,x\,P_m(x)+P_{m-1}(x)
\end{eqnarray}
with initial conditions: $P_{-1}=0,P_0(x)=2.$

{\bf Associated polynomials} (see \ref{associe})

Before computing the associated  polynomials, we need the modified moments, i.e., the moment of the binomial $ 
\bin{x-1}{k}.$

Let us define the linear functional $c^{(s)}$ acting on the space of polynomials as$$\left\langle c^{(s)},x^j\right\rangle :=\frac{B_j}{j!}(s)_{j-1}(-1) ^j,\;\;j\in \N,$$
and $x^2c^{(s)}$ by $$\left\langle x^2c^{(s)},x^j\right\rangle :=\left\langle c^{(s)},x^{j+2}\right\rangle.$$
By recurrence, it is not difficult to prove that
\begin{eqnarray}
\left\langle x^2c^{(2)},
\(\begin{array}{c}
x-1\\
k
\end{array}\)\right\rangle &=&(-1)^k\frac{(k+1)(k+1)!}{(k+3)!}
\end{eqnarray}
and
\begin{eqnarray}
\left\langle x^2c^{(2)},
\(\begin{array}{c}
x-1\\
k
\end{array}\)\(\begin{array}{c}
x-1\\
j
\end{array}\)\right\rangle &=&(-1)^{k+j}\frac{(k+1)(k+1)!(j+1)(j+1)!}{(k+j+3)!}.
\end{eqnarray}

If we define the following polynomial basis $$e_k(x):=\frac{1}{(k+1)(k+1)!}\(\begin{array}{c}
x-1\\
k
\end{array}\),k\in \N,$$

then
$$\left\langle x^2c^{(2)},
e_k\;e_j\right\rangle =\frac{(-1)^{k+j}}{(k+j+3)!}.$$
These moments are related with the coefficients of the exponential function
$$g(x):=\frac{e^{-x}-(1-x+x^2/2)}{x^3}.$$
So we can  recover the expression (\ref{polorths=2}) of the orthogonal polynomials for the functional $x^2c^{(2)}$ by substituting in the orthogonal polynomials for the function $g$
which are 
$$
\sum_{k=0}^{m}\(\begin{array}{c}
m+1\\
k+1
\end{array}\)
\(\begin{array}{c}
m+k+2\\
k+1
\end{array}\)
x^k $$
the monomials $x^k$ by the $e_k's$ (see \ref{hadamard}).

The associated polynomials $R_{m-1}$ (\ref{associe}) are defined as $$R_{m-1}(t) := 
<x^ 2c ^ {(2)},\frac{P_m(x)-P_m(t)}{x-t}>$$
where the variable is $x$.
From the expression (\ref{polorths=2}) for $P_m$, we get the following formula for $R_{m-1}$
\begin{eqnarray}
R_{m-1}=\sum_{k=0}^{m}\(\begin{array}{c}
m+1\\
k+1
\end{array}\)
\(\begin{array}{c}
m+k+2\\
k+1
\end{array}\)\left\langle x^2c^{(2)},\frac{\(\begin{array}{c}
	x-1\\
	k
	\end{array}\)-\(\begin{array}{c}
	t-1\\
	k
	\end{array}\)}{x-t}\right\rangle.
\end{eqnarray}

Using the expression of the polynomial  $\dis\frac{\(\begin{array}{c}
	x-1\\
	k
	\end{array}\)-\(\begin{array}{c}
	t-1\\
	k
	\end{array}\)}{x-t}$
on the Newton basis on $0,1,2,\cdots,k-1$

$$\frac{\bin{x-1}{k}-\bin{t-1}{k}}{x-t}=\bin{t-1}{k}\sum_{i=1} ^k
\frac{\bin{x-1}{i-1}}{i\bin{t-1}{i}},$$
we can write a compact formula for $R_{n-1}$:
$$R_{m-1}(t)=\sum_{k=0}^m
\bin{m+1}{k+1}\bin{m+k+2}{k+1}\bin{t-1}{k}\sum_{i=1}^k
\frac{(-1)^{i-1}}{\bin{t-1}{i}(i+1)(i+2)}
.$$

So, we get the $[m+1/m]$ Pad\'e approximant to the function $\Psi(2,t)$:
$$[m+1/m]_{\psi(2,t)}=B_0-B_1t+t^2\frac{\widetilde{R}_{m-1}(t)}{\widetilde{P}_m(t)}=
B_0-B_1t+t\frac{{R}_{m-1}(1/t)}{{P}_m(1/t)}
$$
and an approximation of $\ze(2,a)$

\begin{eqnarray*}
\ze (2,a) &\approx &\dis{\sum ^{n-1} _{k=0}}
\frac{1}{(k+a)^2}+ \frac{1}{n+a}[m+1/m]_{\Psi(2,t)}\vert_{t=1/(n+a)}  \\&=&
\dis{\sum ^{n-1} _{k=0}}
\frac{1}{(k+a)^2} + \frac{1}{n+a}+\frac{1}{2(n+a)^2}+\frac{1}{(n+a)^2}\varepsilon_m(a)=:\dfrac{v_{n,m}(a)}{u_{n,m}(a)}
\end{eqnarray*}

where 
$$\varepsilon_m(a)=\frac{\sum_{k=0}^m
	\bin{m+1}{k+1}\bin{m+k+2}{k+1}\sum_{i=1}^k
	\frac{\bin{n+a-i-1}{k-i}(-1)^{i-1}}{\bin{k}{i}(i+1)(i+2)}}{\sum_{k=0}^{m}\bin{m+1}{k+1}
	\bin{m+k+2}{k+1}
	\bin{n+a-1}{k}}.$$

{\bf Irrationality of $\ze(2)$}

A consequence of the previous formula is another  proof of the well known irrationality of $\ze(2)$ since it equals to $\pi^2/6.$

Actually, if $d_k:=LCM[1,2,\cdots,k]$, then (\cite{prevost1996}) $$\frac{d_k }{i\bin{k}{i}}\in \N.$$
Then, the numerator of $\varepsilon_m(1)$  multiplied by $d_{m+2}^2$ is an integer and for all integers $n,m$, $d_{m+2}^2v_{n,m}(1)\in \N,d_{m+2}^2u_{n,m}(1) \in \N.$

The error (formula \ref{errorpade}) applied to the function $\Psi(2,t)$ becomes 
$$
\Psi(2,t) - [m+1/m]_{\Psi(2,t)}(t) = \displaystyle{\frac{t^{2m}}{\widetilde
		P^2_m(t)}}\int _{i\R}
{\frac{P^2_m(x)}{1-xt}} \frac{x^2}{\sin^2(\pi x)}dx$$
and the error term satisfies
\begin{eqnarray}
\left| \ze(2)-\frac{v_{n,m}(1)}{u_{n,m}(1)}\right| &\leq& \frac{1}{P^2_m(n+1)}\int _{i\R}
\left| {\frac{P^2_m(x)}{1-x/(n+1)}} \right|  \frac{x^2}{\sin^2(\pi x)} dx\\
&\leq&\frac{1}{P^2_m(n+1)}\int _{i\R}
\left| {P^2_m(x)} \right|  \frac{x^2}{\sin^2(\pi x)} dx\\
&\leq&\frac{1}{P^2_m(n+1)}\frac{2(m+1)(m+2)}{\pi(2m+3)}
\end{eqnarray}

Now,  we consider  $r \in \Q$ such that $m=r n\in\N$. Applying the Stirling formula to the expression  (\ref{polorths=2}) for orthogonal polynomials $P_n$, we get
\begin{eqnarray*} 
  \limsup_n(P_{rn}(n+1)^{1/n}&=&\max_{t\in [0,1]}\frac{(r+t)^{r+t}}{t^{3t}(r-t)^{r- t}(1-t)^{1-t}}\\
  &=&\frac{(r+\sigma(r))^r}{(1-\sigma(r))(r-\sigma(r))^r}=\rho(r)
  \end{eqnarray*}
  where $\sigma(r)=\frac{-r^2+\sqrt{r^4+4r^2}}{2}$ is a zero of $t^2+r^2t-r^2=0$.

  So
  
  \begin{eqnarray*}
   \limsup_n
   \left| 
   \ze(2)
   d_{r n+2}^2u_{n,rn}(1)
   -d_{rn+2}^2v_{n,rn}(1)
   \right|
   ^{1/n} 
   &\le & 
   \limsup_n (d_{rn+2}^2u_{n,rn}(1))^{1/n}
   \limsup_n\frac{1}{(u_{n,rn}(1))^{2/n}}\\
   &=&
  e^{2\max(r,1)}/\rho(r)
    \end{eqnarray*}

The  inequation  $e^{2\max(r,1)}/\rho(r)<1$ is satisfied for $r\in]0.74,1.53[.$

  The irrationality of $\ze(2)$ follows from the 
  following limit
  $$\lim_n\( \ze(2)d_{rn+2}^2u_{n,rn}(1)-d_{rn+2}^2v_{n,rn}(1)\)=0.$$

\subsection{Case $s=3$}

For $s=3$, the weight in the expression  (\ref{zeta}) is
\begin{eqnarray}
w_3(x):=- \frac{x^3}{3}  \int_x^\infty  \frac{d^3}{dt^3}\left(\frac{1}
{e^{2\pi t}-1}\right)dt=\frac{\pi^2 x^3 \cosh \pi x  }{3\sinh^3 (\pi x)}.
\end{eqnarray}
Another expression of $x^2w_3(x)$ is 
$$x^2w_3(x)=\frac{1}{12\pi}\frac{\left| \Gamma(1+i x)\right|^8}{\left| \Gamma(2i x) \right|  ^2}.$$  
This weight has been investigated by Wilson \cite{wilson}.
The orthogonal polynomial $P_n$ satisfies:
\begin{eqnarray}
P_{m}(x^2)&=&(m+1)(m+2)_4F_3\(\begin{array}{cccccc}
-m&m+3&1-x&1+x\\
&&&&;1\\
&2&2&2
\end{array}\)\\&=&
\sum_{k=0}^{m}\bin{m+1}{k+1}
\bin{m+k+2}{k+1}
\bin{x-1}{k}\bin{x+k}{k}/(k+1) \label{polorths=3}
\end{eqnarray}

and
\begin{eqnarray}
\int_{i \R}P_{n}(x^2)P_{m}(x^2)\frac{\pi^2}{3}\frac{x^5\cos \pi x}{\sin^3\pi x}dx&=&0 \;\;(n\neq m)\\
&=&\frac{1}{3}\frac{(-1)^n (n+1)(n+2)}{2n+3}\;\;(n=m)
\end{eqnarray}

The roots of the polynomials $P_n$ are located on the imaginary axis because the weight $\dis\frac{x^5\cos \pi x}{\sin^3\pi x} $is positive on this line.

We set $\Pi_m(x):=P_n(x^2).$

The three terms recurrence relation is
\begin{eqnarray}
\Pi_{m+1}(x)&=&\(\frac{2(2m+3)}{(m+1)(m+2)^2}x^2+\frac{2m+3}{m+2}\)\,\Pi_m(x)-\frac{m+1}{m+2}\Pi_{m-1}(x)
\end{eqnarray}
with initial conditions: $\Pi_{-1}=0,\Pi_0=2.$

{\bf Associated polynomials}

Before computing the associated  polynomials, we need the modified moments, i.e., the moment of the product of  binomials $ 
\bin{x-1}{k}\bin{x+k}{k}.$

In the previous subsection, we have define the linear functional $c^{(s)}$ acting on the space of polynomials as $$\left\langle c^{(s)},x^j\right\rangle :=\frac{B_j}{j!}(s)_{j-1}(-1) ^j,\;\;j\in \N,$$
and $x^4c^{(s)}$ by $$\left\langle x^4c^{(s)},x^j\right\rangle :=\left\langle c^{(s)},x^{j+4}\right\rangle.$$
By recurrence, we can prove 
\begin{eqnarray*}
\left\langle x^4c^{(3)},
\bin{x-1}{k}\bin{x+k}{k}\right\rangle &=&(-1)^{k+1}\frac{(k+1)^2}{2(k+3)(k+2)}
\end{eqnarray*}
and
\begin{eqnarray*}
\left\langle x^4c^{(3)},
\bin{x-1}{k}\bin{x+k}{k}\bin{x-1}{j}\bin{x+j}{j}\right\rangle &=&(-1)^{k+j+1}\frac{(k+1)^2(j+1)^2}{2(k+j+3)(k+j+2)}.
\end{eqnarray*}

If we define the following polynomial basis $$\theta_k(x):=\frac{(-1)^k}{(k+1)^2}\bin{x-1}{k}\bin{x+k}{k}$$

then
$$\left\langle x^4c^{(3)},
			\theta_k\;\theta_j\right\rangle =\frac{(-1)}{2(k+j+3)(k+j+2)}=\frac{1}{2(k+j+3)}-\frac{1}{2(k+j+2)}.$$
These moments are those of the weight function $x(1-x)$ on the interval $[0,1]$.
 
So we can  recover the expression of orthogonal polynomials for the functional $x^4c^{(3)}$ by substituting in the Jacobi orthogonal polynomials  with parameters $\alpha=1,\beta=1$ .
 
$$
 P_m^{(1,1)}(2x-1)=\sum_{k=0}^m\bin{m+1}{k+1}\bin{m+k+2}{k} (-1)^kx^k $$
the monomials $x^k$ by the $\theta_k's$(see \ref{hadamard}).

The associated polynomials $\Theta_{m-1}$ (\ref{associe}) are defined as $$\Theta_{m-1}(t) := 
<x^ 4c ^ {(3)},\frac{\Pi_m(x)-\Pi_m(t)}{x-t}>$$
where the variable is $x$.
From the expression (\ref{polorths=3}) for $\Pi_m$, we get the following formula for $\Theta_{m-1}\in \mathcal{P}_{2m-2}
$
\begin{eqnarray}
\Theta_{m-1}=\sum_{k=0}^{m}\(\begin{array}{c}
m+1\\
k+1
\end{array}\)
\bin{m+k+2}{k+1}
\frac{1}{k+1}
\left\langle x^4c^{(3)},
\frac{
	\bin{x-1}{k}\bin{x+k}{k}-\bin{t-1}{k}\bin{t+k}{k}
	}
	{x-t}
	\right\rangle.
\end{eqnarray}

Using the expression of the polynomial  $\dis
\frac{
	\bin{x-1}{k}\bin{x+k}{k}-\bin{t-1}{k}\bin{t+k}{k}
	}
	{x-t}$
on the Newton basis on $0,1,-1,2,-2\cdots,k,-k$

$$\dis\frac{
	\bin{x-1}{k}\bin{x+k}{k}-\bin{t-1}{k}\bin{t+k}{k}}{x-t}
=
\sum_{j=0}^{2k-1}
\frac
{\nu_{2k}(t)}
{\nu_{j+1}(t)}
\frac
{\nu_j(x)}
{[(j+2)/2]}$$

where
\begin{eqnarray}
\nu_{2j+1}(x)&=&\bin{x-1}{j+1}\bin{x+j}{j}\\
\nu_{2j}(x)&=&\bin{x-1}{j}\bin{x+j}{j}.
\end{eqnarray}

Using
\begin{eqnarray*}
	\left\langle x^4c^{(3)},
\nu_{2j+1}(x)\right\rangle &=-\left\langle x^4c^{(3)},
\nu_{2j}(x)\right\rangle=&(-1)^{j}\frac{(j+1)^2}{2(j+3)(j+2)},
\end{eqnarray*}

we can write a compact formula for $\Theta_{m-1}$:
\begin{eqnarray*}
\Theta_{m-1}(t)&\!\!\!\!\!=&\!\!\!\!\!\sum_{k=0}^m
\bin{m+1}{k+1}\bin{m+k+2}{k+1}\frac{1}{k+1}\sum_{j=0}^{2k-1}
\frac{\nu_{2k}(t)}{\nu_{j+1}(t)}\left\langle x^4c^{(3)},\nu_j(x)\right\rangle\\
=&&\sum_{k=0}^m
\bin{\!m\!+\!1\!}{\!k\!+\!1\!}\bin{\!m\!+\!k\!+\!2\!}{\!k\!+\!1}\frac{\nu_{\!2k\!}(t)}{\!k\!+\!1\!}\sum_{p=1}^{k}
\frac{t}{\bin{\!t\!-\!1\!}{p}\bin{\!t\!+\!p\!}{p}}\frac{(-1)^p}{2(\!p\!+\!1\!)
	(\!p\!+\!2\!)}
\end{eqnarray*}

So, for $m\geq 2$, we get the $[2m+2/2m]$ Pad\'e approximant to the function $\Psi(3,t)$,
$$[2m+2/2m]_{\Psi(3,t)}=\frac{B_0}{2}-B_1t+\frac{3}{2}t^2B_2+t^5\frac{\widetilde{\Theta}_{m-1}(t)}{\widetilde{\Pi}_m(t)}=
\frac{B_0}{2}-B_1t+\frac{3}{2}t^2B_2+t^4\frac{     {\Theta}_{m-1}(1/t)}{ {\Pi}_m(1/t)}
$$
and an approximation of $\ze(3,a)$
$$ \ze (3,a) \approx \dis{\sum ^{n-1} _{k=0}}
\frac{1}{(k+a)^3} + \frac{1}{2(n+a)^2}+\frac{1}{2(n+a)^3}+\frac{1}{4(n+a)^4}+\frac{1}{(n+a)^5}\epsilon_m(a)=:\dfrac{f_{n,m}(a)}{g_{n,m}(a)}
$$

where 
\begin{eqnarray*}
\epsilon_m(a)&=&
	\frac{\Theta_{m-1}(n+a)}
{\Pi_m(n+a)}\\
&=&
\frac{
\sum_{k=0}^m\bin{m+1}{k+1}\bin{m+k+2}{k+1}\frac{n+a}{k+1}
	\sum_{p=1}^k
\frac{	\bin{n+a-p-1}{k-p}\bin{n+a+k}{k-p}(-1)^p}
{\bin{k}{p}^2 2(p+1)(p+2)}
}
{
	\sum_{k=0}^{m}\bin{m+1}{k+1}\bin{m+k+2}{k+1}\bin{n+a-1}{k}\bin{n+a+k}{k}/(k+1)	}.
\end{eqnarray*}

{\bf Irrationality of $\ze(3)$}

The irrationality of $\ze(3)$ has been proved by Apery in a celebrated paper \cite{apery}.
A little later, a particular straightforward and elegant proof of this irrationality was given by Beukers \cite{beukers}.
The author gave another proof in \cite{prevost1996} using the RPA. Actually, $RPA(m,2m-1,2m)$ provides exactly the Apery's numbers.
In the following, we show that the previous approximation also proves irrationality of $\ze(3)$ with RPA's of various degree.

Actually, $$\forall n,a \in \N,d_{m+1}^3\Theta_{m-1}(n+a)\in \Z$$
and $d_{m+1}^3 f_{n,m}(n+a)\in \Z,d_{m+1}^3g_{n,m}(n+a)\in \Z.$

The proof is similar than in previous subsection.

We take $a=1$.

The error (formula \ref{errorpade}) applied to the function $\Psi(3,t)$ becomes 
$$
\Psi(3,t) - [2m+2/2m]_{\Psi(3,t)} = \displaystyle{\frac{t^{2m}}{\widetilde
		\Pi^2_m(t)}}\int _{i\R}
{\frac{\Pi^2_m(x)}{1-xt}} \frac{\pi ^2 x^5\cos\pi x}{3 \sin^3(\pi x)}dx$$
and the error term
$$\ze(3)-\frac{f_{n,m}(1)}{g_{n,m}(1)}$$ satisfies

\begin{eqnarray}
\left| \ze(3)-\frac{f_{n,m}(1)}{g_{n,m}(1)}\right| &\leq&
 \dfrac{1}{\Pi^2_m(n+1)} \int _{i\R}
 \frac{\Pi^2_m(x)}{\left| 1-x/(n+1)\right| } \frac{\pi ^2 x^5\cos\pi x}{3 \sin^3(\pi x)}dx  \\
 &\leq&\dfrac{1}{\Pi^2_m(n+1)} \int _{i\R}
 \Pi^2_m(x) \frac{\pi ^2 x^5\cos\pi x}{3 \sin^3(\pi x)}dx  \\
 &\leq &
\dfrac{1}{\Pi^2_m(n+1)}\frac{1}{3}\frac{ (m+1)(m+2)}{2m+3}
\end{eqnarray}

Now,  we consider  $r \in \Q$ such that $m=r n\in\N$. Using the Stirling formula in the expression of the orthogonal polynomials (\ref{polorths=3}),
\begin{eqnarray*} 
	\limsup_n(\Pi_{rn}(n+1)^{1/n}&=&\max_{t\in [0,1]}\frac{(r+t)^{r+t}(1+t)^{1+t}}{t^{3t}(r-t)^{r- t}(1-t)^{1-t}}\\
	&=&\frac{(1+\mu(r))(r+\mu(r))^r}{(1-\mu(r))(r-\mu(r))^r}=\eta(r)
\end{eqnarray*}
where $\dis\mu(r)=\frac{r}{\sqrt{1+r^2}}$ is a zero of $(1-t^2)(r^2-t^2)=t^4$.

So
\begin{eqnarray*}
	\limsup_n
	\left| 
	\ze(3)
	d_{r n+1}^3 g_{n,rn}(1)
	-d_{rn+1}^3 f_{n,rn}(1)
	\right|
	^{1/n} 
	&\le & 
	\limsup_n (d_{rn+1}^3g_{n,rn}(1))^{1/n}
	\limsup_n\frac{1}{(g_{n,rn}(1))^{2/n}}\\
	&=&
	e^{3\max(r,1)}/\eta(r)
\end{eqnarray*}

The  inequation  $e^{3\max(r,1)}/\eta(r)<1$ is satisfied for $r\in]0.74,1.36[.$

The irrationality of $\ze(3)$ then follows from the 
following limit
$$\lim_n\( \ze(3)d_{rn+1}^3g_{n,rn}(1)-d_{rn+1}^3f_{n,rn}(1)\)=0.$$

 \bibliographystyle{acm}
 \bibliography{bibfile3}

\end{document}